\documentclass[12pt]{article}
\usepackage{amsmath,amsfonts,amssymb,amsthm}
\newtheorem{proposition}{Proposition}[section]
\newtheorem{theorem}[proposition]{Theorem}

\newtheorem{remark}[proposition]{Remark}

\newcommand{\nc}{\newcommand}
\nc{\I}{{\mathbf 1}}

\nc{\bN}{{\mathbf N}}
\nc{\bM}{{\mathbf M}}
\nc{\cB}{{\mathcal B}}
\nc{\cM}{{\mathcal M}}
\nc{\R}{{\mathbb R}}
\nc{\N}{{\mathbb N}}
\nc{\Z}{{\mathbb Z}}
\nc{\BX}{{\mathbb X}}
\nc{\BY}{{\mathbb Y}}
\nc{\cX}{{\mathcal X}}
\nc{\cY}{{\mathcal Y}}

\nc{\BP}{\mathbb{P}}
\nc{\BE}{\mathbb{E}}
\nc{\BQ}{\mathbb{Q}}

\nc{\Po}{{\mathrm{Po}}}
\nc{\Be}{{\mathrm{Beta}}}
\nc{\Di}{{\mathrm{Dir}}}

\numberwithin{equation}{section}

\begin{document} 

\renewcommand{\thefootnote}{\fnsymbol{footnote}}
\author{{\sc G\"unter Last\footnotemark[1]}}
\footnotetext[1]{guenter.last@kit.edu,  
Karlsruhe Institute of Technology, Institute of Stochastics,
76131 Karlsruhe, Germany. }

\title{An integral characterization of the Dirichlet process} 
\date{\today}
\maketitle

\begin{abstract} 
\noindent 
We give a new integral characterization of the Dirichlet process
on a general phase space. To do so we first prove
a characterization of the nonsymmetric Beta distribution
via size-biased sampling. Two applications are a new
characterization of the Dirichlet distribution and a marked version
of a classical characterization of the Poisson-Dirichlet distribution
via invariance under size-biased sampling.
\end{abstract}

\bigskip

\noindent
{\bf Keywords:} 
Dirichlet process, Dirichlet distribution, Beta distribution, Poisson process,
Mecke equation, Poisson-Dirichlet distribution, size-biased sampling

\vspace{0.1cm}
\noindent
{\bf AMS MSC 2010:} 60G55, 60G57

\section{Introduction}\label{intro}

A very important example of a random probability
measure is the Dirichlet process introduced in \cite{Ferguson74}.
More information on this random measure can be found, for instance,
in \cite{Kingman75,Kingman93,Pitman95,Pitman06}.
In this note we study this process on an arbitrary
measurable space $(\BX,\cX)$, equipped with a measure $\rho$ satisfying
$0<\rho(\BX)<\infty$. 

A random measure 
$\zeta$ on $\BX$ (a random element of the space $\bM$ to be defined in Section \ref{sMecke}) 
is called a {\em Dirichlet process} with
{\em parameter measure} $\rho$ if
%if $\BP(\zeta(\BX)=1)=1$, 
%$\BP(\zeta(B)=0)=1$ whenever $\rho(B)=0$ and 
\begin{align}\label{e1.1}
\BP((\zeta(B_1),\ldots,\zeta(B_n))\in\cdot)=\Di(\rho(B_1),\ldots,\rho(B_n)),
\end{align}
whenever $B_1,\ldots,B_n$, $n\ge 1$, form a measurable partition of $\BX$. 
%such that $\rho(B_i)>0$ for each $i\in\{1,\ldots,n\}$. 
We recall that 
the {\em Dirichlet distribution} $\Di(\alpha_1,\dots,\alpha_n)$
with $n\ge 1$ parameters $\alpha_1,\dots,\alpha_n\in[0,\infty)$
such that $\alpha_1+\cdots+\alpha_n>0$ is  
the probability measure on 
$$
\Delta_n:=\{(p_1,\dots,p_n)\in[0,1]^n:p_1+\dots+p_n=1\}
$$
defined as follows. If $n\ge 2$ and $\alpha_1,\dots,\alpha_n\in(0,\infty)$
then $\Di(\alpha_1,\dots,\alpha_n)$ has the density
\begin{align*}
(x_1,\ldots,x_n)\mapsto \frac{\Gamma(\alpha_1+\dots+\alpha_n)}{\Gamma(\alpha_1)\cdots\Gamma(\alpha_n)}
x_1^{\alpha_1-1}\cdots x_{n}^{\alpha_{n}-1}
\end{align*}
with respect to Lebesgue measure on $\Delta_n$.
Here $\Gamma(\alpha):=\int^\infty_0t^{a-1}e^{-t}\,dt$, $\alpha>0$, denotes the Gamma function. 
If $n=1$ then $\Di(\alpha_1):=\delta_1$. If $n\ge 2$ and 
$\alpha_1=\cdots=\alpha_k=0$ for some $k\le n-1$ then
$\Di(\alpha_1,\dots,\alpha_n):=\delta_0^{\otimes k}\otimes \Di(\alpha_{k+1},\dots,\alpha_n)$.
A similar definition applies if some other set of the $\alpha_i$ vanish.
We also recall that the {\em Beta distribution} with parameters $\alpha,\beta>0$
is the probability measure on $[0,1]$ with density
$x\mapsto B(\alpha,\beta)^{-1}x^{\alpha-1}(1-x)^{\beta-1}$,
where 
\begin{align*}
B(\alpha,\beta):=\int^1_0 x^{\alpha-1}(1-x)^{\beta-1}\,dx
=\frac{\Gamma(\alpha)\Gamma(\beta)}{\Gamma(\alpha+\beta)}
\end{align*}
is the Beta function.

Assume that $\zeta$ is a Dirichlet process with parameter measure
$\rho$ and write $\rho=\alpha\nu$, where
$\alpha:=\rho(\BX)$ and $\nu$ is a probability measure on $\BX$.
If $B\in\cX$ has $\nu(B)\in(0,1)$, then $\zeta(B)$ has the Beta
distribution $\Be(\rho(B),\alpha-\rho(B))$.
In particular we have $\BE\zeta(B)=\nu(B)$ for all $B\in\cX$.
It was proved in \cite{Ferguson74} that 
\begin{align}\label{Mecke}
\BE\int f(\zeta,x)\,\zeta(dx)=\BE\iint f((1-u)\zeta+u\delta_x,x)\,G(du)\,\nu(dx),
\end{align}
for all measurable $f\colon\bM\times\BX\to[0,\infty]$,
where $G=\Be(1,\alpha)$.
Under some additional assumptions on $(\BX,\cX,\rho)$ this
was reproved in \cite{SchiLyt17} using the Mecke equation from \cite{Mecke67}.
For the convenience of the reader we provide the proof in Section \ref{sMecke};
see Theorem \ref{Mecke1}.

We say that a probability measure
$\nu$ on $\BX$ is {\em good} if 
%it is either 
%{\em degenerate} (that is $\nu(B)\in\{0,1\}$ for all $B\in\cX$) or 
there exists $B\in\cX$ such that $\nu(B)\in(0,1)\setminus\{1/2\}$.
If $\nu=1/2(\delta_{x}+\delta_{y})$ for distinct $x,y\in\BX$, then $\nu$ is not good.
If $(\BX,\cX)$ is a Borel space these are the only examples
of non-degenerate probability measures which are not good. 
In this note we prove the following characterization of
the Dirichlet process.

\begin{theorem}\label{tMecke} Let $\nu$ be a good probability measure
on $\BX$ and let $G$ be a probability measure on $[0,1]$ satisfying
$b_1:=\int u\,G(du)\in(0,1)$. Assume that $\zeta$ is a random measure
on $\BX$ satisfying \eqref{Mecke} for all
measurable functions $f\colon\bM\times\BX\to[0,\infty)$.
Then $\zeta$ is a Dirichlet process with parameter measure
$\alpha\nu$, where $\alpha:=(1-b_1)/b_1$. Moreover, 
$G=\Be(1,\alpha)$.
\end{theorem}

Theorem \ref{tMecke} has to be compared with 
\cite[Theorem 3.4]{Sethuraman94}, whose proof shows
that a random probability measure $\zeta$ on $\BX$ satisfying
\begin{align}\label{Mecke2}
\BE f(\zeta)%\zeta(\BX)
=\BE\iint f((1-u)\zeta+u\delta_x)\,\Be(1,\alpha)(du)\,\nu(dx),
\end{align}
for all measurable $f\colon\bM\to[0,\infty)$, is a Dirichlet process.
Under some additional assumptions on $(\BX,\cX,\rho)$
this was reproved in \cite{SchiLyt17}, where $\zeta$ was not assumed to
have total mass one. 
Theorem \ref{tMecke} does not make any
assumption on $G$ other than $0<\int u\,G(du)<1$. 
%Moreover, in contrast to much of the
%literature we do not impose any further assumption on the
%probability space $(\BX,\cX,\nu)$.

Let $\zeta$ be a random probability measure satisfying \eqref{Mecke}
and let $\tau$ be a random element of $\BX$ satisfying
$\BP(\tau\in B \mid \zeta)=\zeta(B)$ almost surely for each $B\in\cX$.
Then $\tau$ has distribution $\nu$ and %\eqref{Mecke} implies that 
\begin{align}\label{e1.11}
\BP[\zeta\in\cdot\mid \tau=x]=\BP((1-W)\zeta+W\delta_x\in\cdot),\quad \nu\text{-a.e.\ $x\in\BX$},
\end{align}
where  $W$ is independent of $(\zeta,\tau)$ and has distribution $G$.
If $\zeta$ is a Dirichlet process with parameter measure $\rho$, 
then $W$ has distribution $\Be(1,\rho(\BX))$ and
$(1-W)\zeta+W\delta_x$ is a Dirichlet
process with parameter measure $\rho+\delta_x$.
Equation  \eqref{Mecke} shows in particular that
$\zeta\overset{d}{=}(1-W)\zeta+W\delta_X$, where $X$ is independent of 
$(\zeta,W)$ and has distribution $\nu$.
Iterating this identity yields (in the case $G(\{0\})<1$)
\begin{align}\label{e1.12}
\zeta\overset{d}{=}\sum^\infty_{n=1}W_n\prod^{n-1}_{i=1}(1-W_i)\delta_{X_n},
\end{align}
where the sequence $(W_n)$ is i.i.d.\ $G$,
the sequence $(X_n)$ is i.i.d.\ $\nu$, and both sequences are
independent. If $(\BX,\cX)$ is a Borel space and $\nu$ is diffuse,
we may use a classical result from \cite{McCloskey65} 
on {\em size-biased permutations} (see \cite[Theorem 1]{Hoppe86} and
\cite[Theorem 5]{Pitman95})
to obtain from \eqref{e1.11} and \eqref{e1.12} that $G=\Be(1,\alpha)$
for some $\alpha>0$; see also Section \ref{sinv}.
If $\nu$ is not diffuse, we have not been
able to derive Theorem \ref{tMecke} from the literature, even not in the special case of 
a space $\BX$ consisting of just two points. The latter case is treated in
our Theorem \ref{tbeta}, a characterization of the nonsymmetric Beta distribution.
In Section \ref{sDdistr} we specialize Theorem \ref{tMecke} to the
case of the Dirichlet distribution. In Section \ref{sinv} we use
Theorems \ref{tMecke} and its converse (Theorem \ref{Mecke1}) 
to derive a marked version of a classical
characterization of the {\em Poisson-Dirichlet distribution} via invariance
under size-biased sampling. 

Our approach is based on the classical Poisson representation
of a Dirichlet process (recalled in Section \ref{sMecke}) and a few 
fundamental properties of the Poisson process and the Dirichlet
distribution.

\section{The Mecke equation for the Dirichlet process}\label{sMecke}

We start by briefly recalling the classical construction of the Dirichlet process 
in terms of a suitable Poisson process; see \cite{Ferguson74,Kingman75}.
Let $\bM$ denote the space of all $s$-finite measures on $(\BX,\cX)$,
equipped with the smallest $\sigma$-field making the mappings
$\mu\mapsto \mu(B)$ for each $B\in\cX$ measurable on $\bM$. 
A {\em random measure} $\xi$ on $\BX$ is a random element of
$\bM$ defined over a suitable probability space
$(\Omega,\mathcal{F},\BP)$ with associated integral operator $\BE$;
see e.g.\ \cite[Chapter 13]{LastPenrose17}.
Let $\eta$ be a Poisson process on $\BX\times(0,\infty)$ with
intensity measure $\BE\eta(d(x,r))=r^{-1}e^{-r}dr\rho(dx)$ and define
a random measure $\xi$ on $\BX$ by
\begin{align}\label{e14.1}
\xi(B)=\int_{B\times(0,\infty)} r\,\eta(d(y,r)),\quad B\in \cX.
\end{align}
Then $\xi$ is completely independent, and $\xi(B)$ has a Gamma distribution with 
shape parameter $\rho(B)$ and scale parameter $1$; see also
\cite[Example 15.6]{LastPenrose17}. Define a random measure $\zeta$ on $\BX$ by
$\zeta(B):=\xi(B)/\xi(\BX)$, $B\in\cX$, where $0/0:=0$.
Then $\zeta$ is a Dirichlet process with parameter measure $\rho$, independent of $\xi(\BX)$.
By \cite[Proposition 13.2]{LastPenrose17} the distribution of a Dirichlet process
is determined by its parameter measure. 
%Hence we can assume without loss of generality
%that it is defined in terms of $\eta$.

The next result was proved in \cite{Ferguson74} using a different formulation.
Here we follow \cite{SchiLyt17} in providing a short proof,
based on the classical Mecke equation from \cite{Mecke67}.

\begin{theorem}\label{Mecke1} Let $\zeta$ be a Dirichlet
process on $\BX$ with parameter measure $\rho$ and let $f\colon\bM\times\BX\to[0,\infty)$ be measurable.
Then \eqref{Mecke} holds with $G=\Be(1,\alpha)$, where $\alpha:=\rho(\BX)$ and
$\nu:=\alpha^{-1}\rho$.
\end{theorem}
{\em Proof.} Since \eqref{Mecke} does only concern the distribution
of $\zeta$, we can assume that $\zeta$ is defined in terms of
the Poisson process $\eta$ as above.
Since $\BE\xi(\BX)=\alpha$ and since
$\zeta$ and $\xi(\BX)$ are independent, we have
\begin{align*}
\BE\int f(\zeta,x)\,\zeta(dx)&=\alpha^{-1}\BE\int f(\zeta,x)\xi(\BX)\,\zeta(dx)\\
&=\alpha^{-1}\BE\int f(\zeta,x)\,\xi(dx)\\
&=\alpha^{-1}\BE\int f(\xi(\BX)^{-1}\xi,x)r\,\eta(d(x,r)).
\end{align*}
Since $\xi$ is a measurable function of the Poisson process $\eta$,
we can apply the Mecke
equation (see \cite{Mecke67} and \cite[Theorem 4.1]{LastPenrose17})
to obtain that
\begin{align*}
\BE\int f(\zeta,x)\,\zeta(dx)
=\alpha^{-1}\BE\iint f\bigg(\frac{\xi+r\delta_x}{\xi(\BX)+r},x\bigg)e^{-r}\,dr\,\rho(dx)\\
=\BE\iiint f\bigg(\frac{s\zeta+r\delta_x}{s+r},x\bigg)e^{-r}\,G(\alpha,1)(ds)\,dr\,\nu(dx),
\end{align*}
where $G(\alpha,1)$ denotes the Gamma distribution with shape parameter $\alpha$ and
scale parameter $1$.
Since $\iint \I\{r/(r+s)\in\cdot\}\,e^{-r}\,G(\alpha,1)(ds)\,dr=\Be(1,\alpha)$,
we obtain the assertion.\qed

%\bigskip

%{\em Proof.} Since \eqref{Mecke} does only concern the distribution
%of $\zeta$, we can assume that $\zeta$ is defined in terms of
%the Poisson process $\eta$ as above.
%Since $\BE\xi(\BX)^n=\alpha_{(n)}:=\alpha\cdots (\alpha+n-1)$ and since
%$\zeta$ and $\xi(\BX)$ are independent, we have
%\begin{align*}
%\BE\int &f(\zeta,x_1,\ldots,x_n)\,\zeta^n(d(x_1,\ldots,x_n))\\
%&=\alpha_{(n)}^{-1}\,\BE\int f(\zeta,x_1,\ldots,x_n)\xi(\BX)^n\,\zeta^n(d(x_1,\ldots,x_n))\\
%&=\alpha^{-1}_{(n)}\,\BE\int f(\zeta,x_1,\ldots,x_n)\,\xi^n(d(x_1,\ldots,x_n))\\
%&=\alpha^{-1}_{(n)}\,\BE\int f(\xi(\BX)^{-1}\xi,r_1x_1,\ldots,r_nx_n)r_1\cdots r_n\,
%\eta^n(d((x_1,r_1),\ldots,(x_n,r_n)).
%\end{align*}
%Since $\xi$ is a measurable function of the Poisson process $\eta$,
%we can apply the multivariate Mecke
%equation (see \cite{Mecke67} and \cite[Theorem 4.4]{LastPenrose17})
%to obtain that
%% the following argument is wrong, use instead that \eta is simple 
%\begin{align*}
%\BE\int f(\zeta,x)\,\zeta(dx)
%=\alpha^{-1}\BE\iint f\bigg(\frac{\xi+r\delta_x}{\xi(\BX)+r},x\bigg)e^{-r}\,dr\,\rho(dx)\\
%=\BE\iiint f\bigg(\frac{s\zeta+r\delta_x}{s+r},x\bigg)e^{-r}\,G(\alpha,1)(ds)\,dr\,\nu(dx).
%\end{align*}

\section{A characterization of the Beta distribution}\label{schbeta}

In this section we prove the following characterization of the nonsymmetric Beta  distribution. 
%given parameters $\alpha_1,\alpha_2>0$. It is convenient to
%rewrite these parameters as  $\alpha:=\alpha_1+\alpha_2$
%and $p:=\alpha_1/\alpha$.

\begin{theorem}\label{tbeta} Let $Z$ be a $[0,1]$-valued random variable
such that $\BE Z\ne 1/2$. Then $Z$ has a Beta distribution
iff there exists an independent $[0,1]$-valued random variable $W$
%a probability measure $G$ on $[0,1]$
with $\BE W\in(0,1)$ and a number $p\in(0,1)\setminus\{1/2\}$ such that
\begin{align}\label{e3.1}
\BE g(Z)Z&=p\,\BE g((1-W)Z+W),\\ \label{e3.2}
\BE g(Z)(1-Z)&=(1-p)\,\BE g((1-W)Z)
\end{align}
for all measurable $g\colon[0,1]\to [0,\infty)$.
In this case $Z$ has the
$\Be(p\alpha,(1-p)\alpha)$ distribution
and $W$ has the $\Be(1,\alpha)$ distribution, where $\alpha=(1-\BE W)/\BE W$. 
\end{theorem}
{\em Proof.} Let $p\in(0,1)$ and $\alpha>0$. Suppose that
$\zeta$ is a Dirichlet process on $\{1,2\}$
with parameter measure $\rho=\alpha((1-p)\delta_1+p\delta_2)$.
Define $Z':=\zeta(\{2\})$ and $G':=\Be(1,\alpha)$.
Applying \eqref{Mecke} with $f(\zeta,x)=g(\zeta(\{2\}))\I\{x=2\}$
yields
\begin{align}\label{e3.3}
\BE g(Z')Z'&=p\,\BE \int g((1-u)Z'+u)\,G'(du),
\end{align}
while the choice $f(\zeta,x)=g(\zeta(\{2\}))\I\{x=1\}$ yields
\begin{align}\label{e3.4}
\BE g(Z')(1-Z')&=(1-p)\,\BE \int g((1-u)Z')\,G'(du).
\end{align}
Since $Z'$ has the $\Be(p\alpha,(1-p)\alpha)$ distribution,
this proves one implication.

Assume now that $W$ and $p\in(0,1)\setminus\{1/2\}$ have the stated properties and define
$\alpha:=(1-\BE W)/\BE W$. Let $Z'$ and $G'$ be given as in the first
part of the proof. We shall show that $Z$ and $Z'$ have the same
distribution and that $W$ has distribution $G'$. 
To this end we prove equality of moments
by induction.
For each $n\in\N$ we define $a_n:=\BE Z^n$,
$a'_n:=\BE Z'^n$, $b_n:=\BE W^n$ and $b'_n:=\int u^n\,G'(du)$.
Choosing $g\equiv 1$ in \eqref{e3.1} gives
$a_1=p$, which equals $a'_1$. Further
$b'_1=1/(\alpha+1)=b_1$. Given $n\in\N$ we assume now that $a_k=a'_k$ and
$b_k=b'_k$ for $k\in\{1,\ldots,n\}$. Choosing $g(x)=x^{n}$ in
\eqref{e3.1} and \eqref{e3.3} shows that
\begin{align}\label{e3.7}
a_{n+1}=a'_{n+1}.
\end{align}
%Assume first that $n$ is odd.
Taking $g(x)=x^{n+1}$ in \eqref{e3.1} and \eqref{e3.2}
and combining both equations yields
\begin{align*}
%a_{n+1}-p\,\BE \int (Z+u(1-Z))^{n+1}\,G(du)=(1-p)\,\BE \int ((1-u)Z)^{n+1}\,G(du),
a_{n+1}-p\,\BE  (Z+W(1-Z))^{n+1}=(1-p)\,\BE ((1-W)Z)^{n+1},
\end{align*}
or, more explicitly, 
\begin{align*}
a_{n+1}-p\sum^{n+1}_{k=0}\binom{n+1}{k} b_{n+1-k}\BE Z^k(1-Z)^{n+1-k}
-(1-p)a_{n+1} \BE (1-W)^{n+1}=0.
\end{align*}
This can be solved for $b_{n+1}$ (in terms of
$a_1,\ldots,a_{n+1},b_1,\ldots,b_n$) if 
\begin{align}\label{e3.9}
p\,\BE (1-Z)^{n+1}+(1-p)a_{n+1}(-1)^{n+1}\ne 0.
\end{align}
By \eqref{e3.3}, \eqref{e3.4}, \eqref{e3.7}
and induction hypothesis we have 
$$
\BE (1-Z)^{n+1}=\BE (1-Z')^{n+1}=\prod^n_{j=0}\frac{p\alpha+j}{\alpha+j},
$$
where the second identity is well-known (and follows from an easy
calculation). Similarly,
$$
a_{n+1}=\BE (Z')^{n+1}=\prod^n_{j=0}\frac{(1-p)\alpha+j}{\alpha+j}.
$$
Therefore \eqref{e3.9} holds if
\begin{align}\label{e3.11} 
p\prod^n_{j=0}(p\alpha+j)+(-1)^{n+1}(1-p)\prod^n_{j=0}((1-p)\alpha+j)\ne 0.
\end{align}
If $n$ is odd, this is true. If $n$ is even, this follows from
$p\ne 1/2$. 
In view of \eqref{e3.3} and \eqref{e3.4}, $\BE (Z')^{n+1}$ can be expressed
in terms of $a_1,\ldots,a_{n+1},b_1,\ldots,b_n$ by the same function
as $b_{n+1}$. This finishes the induction step and hence concludes the
proof.\qed

\bigskip 

Theorem \ref{tbeta} has some similarities with Lemma 13 in \cite{Pitman96},
which is another characterization of the nonsymmetric Beta distribution.

\begin{remark}\rm
At the moment we do not know whether Theorem \ref{tbeta} remains true in
the symmetric case $p=1/2$. In this case it is easy to show
that $Z$ and $1-Z$ have the same distribution. Assuming this
distributional equality, equation \eqref{e3.1} holds for all functions $g$
iff the same is true for \eqref{e3.2}. In essence this raises the following problem.
Assume that $Z\overset{d}{=}1-Z$ and that
$V$ is an independent $[0,1]$-valued random variable
with $\BE V\in(0,1)$ satisfying
\begin{align*}
\BE (1-Z)^n Z=\tfrac{1}{2}\,\BE(1-Z)^n\,\BE V^n,\quad n\in\N_0. 
\end{align*}
Is it then true that $Z$ and $V$ have a Beta distribution? 
\end{remark}

Replacing \eqref{e3.2} by another
condition yields a characterization of the general Beta distribution,
similar to the one given in \cite{Pitman96,SeshWes03}.

\begin{theorem}\label{tbeta2} Let $Z$ be a $[0,1]$-valued random variable.
Then $Z$ has a Beta distribution
iff there exists an independent $[0,1]$-valued random variable $W$
with $\BE W\in(0,1)$ and numbers $p\in(0,1)$ and $c\ge 0$ such that
\eqref{e3.1} and
\begin{align}\label{e3.8}
\BE g(Z)Z^2=c\,\BE g((1-W)Z+W)W
\end{align}
hold for all measurable $g\colon[0,1]\to [0,\infty)$, where
$\alpha:=(1-\BE W)/\BE W$.
In this case $Z$ has the
$\Be(p\alpha,(1-p)\alpha)$ distribution,
$W$ has the $\Be(1,\alpha)$ distribution and $c=p(\alpha p+1)$.
\end{theorem}
{\em Proof.} %We only need to prove one implication and 
Assume
the stated conditions. Choosing $g\equiv 1$ in \eqref{e3.8} yields
$\BE Z^2=c\,\BE W=c/(\alpha+1)$. On the other hand we have from 
\eqref{e3.1} that
\begin{align*}
\BE Z^2=p\,\BE (W+Z-WZ)=p\Big(\frac{1}{\alpha+1}+p-\frac{1}{\alpha+1}\Big)
=p\frac{\alpha p+1}{\alpha+1},
\end{align*}
so that $c=p(\alpha p+1)$.

If \eqref{e3.1} holds, then assumption \eqref{e3.8} is equivalent to
\begin{align*}
p\,\BE g(Z+W-WZ)(Z+W-WZ)=c\,\BE g(Z+W-WZ)W,
\end{align*}
or
\begin{align}\label{e3.41}
\BE g(Z+W-WZ)(Z+W-WZ)=(\alpha p+1)\,\BE g(Z+W-WZ)W.
\end{align}
If $Z$ and $W$ have the stated Beta distributions, this is true.
In fact, $W/(Z+W-WZ)$ and $Z+W-WZ$ are independent in this case; see e.g.\ 
\cite{SeshWes03}.

From now on we assume that \eqref{e3.1} and \eqref{e3.41} hold.
We are using \eqref{e3.41} in the form
\begin{align}\label{e3.42}
\BE g(Z+W-WZ)(Z-WZ)=\alpha p\,\BE g(Z+W-WZ)W.
\end{align}
As in the proof of Theorem \ref{tbeta} we use induction and can
proceed until \eqref{e3.7} without any change.
But now we take $g(x)=x^n$ in \eqref{e3.42} to obtain that
\begin{align*}
\BE (W(1-Z)+Z)^n(1-W)Z=\alpha p\,\BE (W(1-Z)+Z)^n W,
\end{align*}
that is 
\begin{align*}
\sum^n_{k=0}\binom{n}{k}b_k\, \BE (1-Z)^kZ^{n-k+1}-\sum^n_{k=0}\binom{n}{k}b_{k+1}\, &\BE (1-Z)^kZ^{n-k}\\
&=\alpha p\,\sum^n_{k=0}\binom{n}{k}b_{k+1}\BE (1-Z)^kZ^{n-k}.
\end{align*}
This can be uniquely solved for $b_{n+1}$ in terms of $(a_1,\ldots,a_{n+1},b_1,\ldots,b_n)$.\qed

\bigskip

Since \eqref{e3.8} does not follow from \eqref{e3.1} and \eqref{e3.2},
we cannot use Theorem \ref{tbeta2} to prove Theorem \ref{tbeta}
in the symmetric case.

\section{Proof of Theorem \ref{tMecke}}

In this section we prove Theorem \ref{tMecke}, assuming that its assumptions
are satisfied.

Choosing $f\equiv 1$ in \eqref{Mecke} gives $\BE\zeta(\BX)=1$.
Choosing $f(\zeta,x)=\zeta(\BX)$ yields
\begin{align*}
\BE\zeta(\BX)^2=\BE\int ((1-u)\zeta(\BX)+u)\,G(du)
=\int ((1-u)+u)\,G(du)=1.
\end{align*}
Hence $\zeta(\BX)$ has variance $0$, so that
$\BP(\zeta(\BX)=1)=1$. 

%Assume first that $\nu$ is degenerate, that is $\nu(B)=\BE\zeta(B)\in\{0,1\}$
%for all $B\in\cX$. Then, for each $B\in\cX$, either $\BP(\zeta(B)=0)=1$
%or $\BP(\zeta(\BX\setminus B)=0)=1$. Hence, for trivial reasons, $\zeta$ is a
%(degenerate) Dirichlet process with parameter measure $\alpha\nu$.

%Assume now that $\nu$ is not degenerate. 
Since $\nu$ is assumed to be good,
there exists $B\in\cX$ such that $p:=\nu(B)\in(0,1)\setminus\{1/2\}$.
Let $Z:=\zeta(B)$. Choosing $f(\zeta,x)=g(\zeta(B))\I\{x\in B\}$
in \eqref{Mecke} gives \eqref{e3.1}, while
choosing $f(\zeta,x)=g(\zeta(B))\I\{x\notin B\}$ gives \eqref{e3.2}.
Theorem \ref{tbeta} shows that $G=\Be(1,\alpha)$, where $\alpha:=(1-b_1)/b_1$.

Let $n\ge 2$ and let $B_1,\ldots,B_n$ be a measurable partition of $\BX$ 
such that $\rho(B_i)>0$ for each $i\in\{1,\ldots,n\}$.
For $k_1,\ldots,k_n\in\N_0$ we define
\begin{align*}
b(k_1,\ldots,k_n):=\BE \zeta(B_1)^{k_1}\cdots \zeta(B_n)^{k_n}.
\end{align*}
Defining $\alpha_i:=\alpha\nu(B_i)$, $i\in\{1,\ldots,n\}$, we need to show that
\begin{align*}
b(k_1,\ldots,k_n)=b'(k_1,\ldots,k_n):=
\int x^{k_1}_1\cdots x^{k_n}_n\, \Di(\alpha_1,\dots,\alpha_n)(d(x_1,\ldots,x_n)).
\end{align*}
It is well-known (and easy to derive from the Beta integral) that
\begin{align}\label{e411}
b'(k_1,\ldots,k_n)=\frac{\Gamma(\alpha)}{\Gamma(\alpha+k_1+\cdots+k_n)}
\prod^n_{i=1}\frac{\Gamma(\alpha_i+k_i)}{\Gamma(\alpha_i)}.
\end{align}

Let $j\in\{1,\ldots,n\}$.
Choosing 
$$
f(\zeta,x)=\zeta(B_1)^{k_1}\cdots\zeta(B_n)^{k_n}\I\{x\in B_j\}
$$
in \eqref{Mecke} gives us
\begin{align}\label{e4.1}\notag
&b(k_1,\ldots,k_j+1,\ldots,k_n)\\ \notag
&=\nu(B_j)\sum^{k_j}_{r=0}\binom{k_j}{r} b(k_1,\ldots,r,\ldots,k_n)
\int (1-u)^{k_1}\cdots(1-u)^ru^{k_j-r}\cdots (1-u)^{k_n}\,G(du)\\
&=\alpha_j\sum^{k_j}_{r=0}\binom{k_j}{r} b(k_1,\ldots,r,\ldots,k_n)
B(k_j+1-r,k_1+\cdots+k_n+\alpha+r-k_j),
\end{align}
where we have used that $G(du)=\alpha (1-u)^{\alpha-1}du$.
These recursions determine the moments $b(k_1,\ldots,k_n)$.
It remains to show that the moments in
\eqref{e411} do also satisfy these recursions.
We do this for $j=n$ and fixed $k_1,\ldots,k_{n-1}\in\N_0$.
Writing $b_k:=b'(k_1,\ldots,k_{n-1},k)$
and $k_0:=k_1+\cdots +k_{n-1}$, we need to show that
\begin{align*}
b_{k+1}=\alpha_n\sum^{k}_{r=0}\binom{k}{r} b_r
\frac{\Gamma(k+1-r)\Gamma(k_0+\alpha+r)}{\Gamma(\alpha+k_0+k+1)}.
\end{align*}
Inserting here \eqref{e411} and $\Gamma(k+1-r)=(k-r)!$ this boils down
to
\begin{align*}
\frac{\Gamma(\alpha_n+k+1)}{k!}=\alpha_n\sum^{k}_{r=0}\frac{\Gamma(\alpha_n+r)}{\Gamma(r+1)},  
\end{align*}
which can be confirmed by induction.\qed

\bigskip

Assume that $\zeta$ is a random measure on $\BX$ satisfying the assumptions of
Theorem \ref{tMecke}, except that $\nu$ is {\em degenerate}, that is $\nu(B)=\BE\zeta(B)\in\{0,1\}$
for all $B\in\cX$. Then, for each $B\in\cX$, either $\BP(\zeta(B)=0)=1$
or $\BP(\zeta(\BX\setminus B)=0)=1$. Hence, for trivial reasons, 
$\zeta$ is a (degenerate) Dirichlet process with parameter measure $\alpha\nu$.

\section{A characterization of the Dirichlet distribution}\label{sDdistr}

Given $n\in\N$ we denote by $e_1,\ldots,e_n$
the standard basis in $\R^n$.

\begin{theorem}\label{tDirichletc} Let $Z=(Z_1,\ldots,Z_n)$ be a
random probability vector in $\R^n$ whose expectation vector
$(p_1,\ldots,p_n):=\BE Z$ forms a good probability measure on $\{1,\ldots,n\}$.
Then $Z$ has a Dirichlet distribution
iff there exists an independent $[0,1]$-valued random variable $W$
with $\BE W\in(0,1)$ such that
\begin{align}\label{e3.71}
\BE g(Z)Z_i=p_i\,\BE g((1-W)Z+We_i),\quad i=1,\ldots,n,
\end{align}
for all measurable $g\colon[0,1]\to [0,\infty)$. In this case
$W$ has the $\Be(1,\alpha)$ distribution and
$Z$ has the $\Di(\alpha p_1,\ldots,\alpha p_n)$ distribution,
where $\alpha=(1-\BE W)/\BE W$.
\end{theorem}
{\em Proof.} The result follows from Theorems \ref{tMecke}
and \ref{Mecke1} applied to the case $\BX=\{1,\ldots,n\}$
and $\nu=p_1\delta_1+\cdots+p_n\delta_n$.\qed

\bigskip 

In the case $n=2$, Theorem \ref{tDirichletc} boils down to
Theorem \ref{tbeta2}. 

Equation \eqref{e3.71} can be expressed
in different ways. For instance we can consider, as in
\eqref{e1.11}, a $\{1,\ldots,n\}$-valued random variable $\tau$
with $\BP(\tau=i\mid Z)=Z_i$, $i\in \{1,\ldots,n\}$.
Then the equations \eqref{e3.71} hold iff
\begin{align}
\BE f(Z,\tau)=\sum^n_{i=1}p_i\,\BE f((1-W)Z+We_i,i),
\end{align}
for all measurable $f\colon[0,1]\times\{1,\ldots,n\}\to [0,\infty)$.

\section{Invariance under size-biased sampling}\label{sinv}

If $(\BX,\cX)$ is a Borel space and $\mu$ is a measure
on $\BX$ we define
$$
\mu^{(x)}:=(1-\mu\{x\})^{-1}(\mu-\mu\{x\}\delta_x)
$$ 
where $\mu\{x\}:=\mu(\{x\})$ and $a/0:=0$ for all $a\in\R$.
Effectively we will use this definition only for probability
measures $\mu$ with $\mu\{x\}<1$ for all $x\in\BX$.
The following result is a straightforward consequence of 
Theorem \ref{Mecke1}.

\begin{theorem}\label{t5.1} Assume that $(\BX,\cX)$ is a Borel space
and that $\rho$ is a diffuse measure on $\BX$.
Let $\zeta$ be a Dirichlet process with parameter measure $\rho$ and
let $\tau$ be a random element of $\BX$ satisfying 
$\BP(\tau\in B \mid \zeta)=\zeta(B)$ almost surely for each $B\in\cX$.
Then $\zeta^{(\tau)}\overset{d}{=}\zeta$. Moreover $\zeta^{(\tau)}$,
$\zeta\{\tau\}$ and $\tau$ are independent.
\end{theorem}

%For later reference we mention yet another construction of a Dirichlet
%process, based on a Poisson process 
To discuss and to reverse this theorem  we need some notation.
If $x=(x_n)_{n\ge 1}$ is a sequence in some space and $i\in\N$
then we let $x^i=(x^i_n)_{n\ge 1}$ denote the sequence arising from
$x$ by dropping the $i$-th member, that is
$x^i_n:=x_n$ for $n<i$ and $x^i_n:=x_{n+1}$ for $n\ge i$.
If $x$ is real-valued we define $x^{(i)}:=(1-x_i)^{-1}x^{i}$,
where $a/0:=0$ for each $a\in\R$.
%If $x$ is a probability distribution, then so is $x^{(i)}$.

Let $\alpha>0$ and let $\eta'=\sum^\infty_{n=1}\delta_{T_n}$ be 
a Poisson process on $(0,\infty)$ with intensity measure $\BE\eta'(dr)=\alpha r^{-1}e^{-r}dr$.
Let $\zeta':=S^{-1}\sum^\infty_{n=1}\delta_{T_n}$, where $S:=\sum^\infty_{n=1}T_n$
and let $X=(X_n)_{n\ge 1}$ be an independent i.i.d.\ sequence with
marginal distribution $\nu$. The construction in Section \ref{sMecke} 
and the marking theorem
for Poisson processes (see e.g.\ \cite[Theorem 5.6]{LastPenrose17}) show that
$S^{-1}\sum^\infty_{n=1}\delta_{(T_n,X_n)}$ is a Dirichlet process
with parameter measure $\alpha\nu$; see \cite{Ferguson74} for a slightly
different argument. 
By writing $\zeta'=\sum^\infty_{n=1}\delta_{Z_n}$, where $Z_1>Z_2>\cdots$, we can
identify $\zeta'$ with the sequence $Z:=(Z_n)_{n\ge 1}$. 
Both, the random sequence $Z$ and
the random measure $\zeta'$ are said to have
the {\em Poisson-Dirichlet distribution} with {\em parameter} $\alpha$; see \cite{Pitman95}.
Let $\tau$ be a random variable as in Theorem \ref{t5.1} and define
a random element $\kappa$ of $\N$ by $\kappa:=i$ if $\tau=X_i$. Then
\begin{align}\label{kappa}
\BP(\kappa=i\mid Z,X)=Z_i,\quad \text{$\BP$-a.s.},\,i\in\N.
\end{align}
Moreover, $\zeta^{(\tau)}$ can be identified with $(Z^{(\kappa)},X^{\kappa})$,
while $\zeta\{\tau\}=Z_\kappa$ and $\tau=X_\kappa$.
Hence Theorem \ref{t5.1} shows that
$(Z^{(\kappa)},X^{\kappa})\overset{d}{=}(Z,X)$ and that
$(Z^{(\kappa)},X^{\kappa})$, $Z_\kappa$ and $X_\kappa$ are independent.
This can be reversed without assuming that $\rho$ is diffuse.

\begin{theorem}\label{t5.2} Suppose that
$Z=(Z_n)_{n\ge 1}$ is a $(0,1)$-valued sequence with $\sum^\infty_{n=1} Z_n=1$
and that $X=(X_n)_{n\ge 1}$ is a sequence of random elements of $\BX$.
Let $\kappa$ be a random element of $\N$ satisfying \eqref{kappa}.
Assume that $(Z^{(\kappa)},X^{\kappa})\overset{d}{=}(Z,X)$ and that
$(Z^{(\kappa)},X^{\kappa})$, $Z_\kappa$ and $X_\kappa$ are independent.
Then $Z_\kappa$ has the $\Be(1,\alpha)$ distribution, where 
$\alpha:=(1-\BE Z_\kappa)/\BE Z_\kappa$. Moreover,
$\sum^\infty_{n=1}Z_n\delta_{X_n}$ is a Dirichlet process
with parameter measure $\alpha\nu'$, where $\nu':=\BP(X_\kappa\in\cdot)$.
If $\BX$ is a Borel space, then $\sum^\infty_{n=1}\delta_{Z_n}$ has the
Poisson-Dirichlet distribution with parameter $\alpha$.
\end{theorem}
{\em Proof.} Let $\lambda$ be a diffuse probability measure on some 
Borel space $(\BY,\cY)$ and let $Y=(Y_n)_{n\ge 1}$ be an i.i.d.\ sequence with
marginal distribution $\lambda$, independent
of $(Z,X,\kappa)$. Then
\begin{align*}
\zeta:=\sum^\infty_{n=1}Z_n\delta_{(X_n,Y_n)}
\end{align*}
is a random probability measure on $\BX\times\BY$.
Note that $\zeta=T(Z,X,Y)$ for a well-defined measurable mapping $T$. 
Our goal is to establish the integral equation \eqref{Mecke}
with $\BX$ replaced by $\BX\times\BY$, $\nu:=\nu'\otimes \lambda$
and $G:=\BP(Z_\kappa\in\cdot)$. 

Taking a measurable function $f$ with suitable domain we have
\begin{align}\label{e5.3}
\BE\int f(\zeta,v,w)\,\zeta(d(v,w))=
\sum^\infty_{i=1}\BE \int f(T(Z,X,y),X_i,y_i)Z_i\,\lambda^\infty(dy).
\end{align}
For each $i\in\N$ and each sequence $y=(y_n)_{n\in\N}\in\BY^\infty$ we have
\begin{align*}
T(Z,X,y)=\sum_{n\ne i}Z_n\delta_{(X_n,y_n)}+ Z_i\delta_{(X_i,y_i)}
=(1-Z_i)T(Z^{(i)},X^i,y^i)+Z_i\delta_{(X_i,y_i)},
\end{align*}
so that the right-hand side of \eqref{e5.3} equals
\begin{align*}
\sum^\infty_{i=1}\BE& \iint f((1-Z_i)T(Z^{(i)},X^i,y),X_i,w)+Z_i\delta_{(X_i,w)},X_i,w)
Z_i\,\lambda^\infty(dy)\,\lambda(dw)\\
&=\iint \BE f((1-Z_\kappa)T(Z^{(\kappa)},X^{\kappa},y),X_\kappa,w)+Z_\kappa\delta_{(X_\kappa,w)},X_\kappa,w)
\,\lambda^\infty(dy)\,\lambda(dw).
\end{align*}
By our assumptions this equals
\begin{align*}
\iiiint \BE f((1-r)T((Z,X,y),v,w)+r\delta_{(v,w)},v,w)\,G(dr)\,\nu'(dv)\,
\,\lambda^\infty(dy)\,\lambda(dw).
\end{align*}
By the definition of $\zeta$, this implies the desired
integral equation \eqref{Mecke}.

Since $\nu'\otimes\lambda$ is good, Theorem \ref{tMecke} shows that $G=\Be(1,\alpha)$
and that $\zeta$ is a Dirichlet process
with parameter measure $\alpha\nu'\otimes\lambda$. 
In particular, $\sum^\infty_{n=1}Z_n\delta_{X_n}$ is a Dirichlet process
with parameter measure $\alpha\nu'$. Moreover, by the construction
given above, 
\begin{align}\label{e5.22}
\zeta=\sum^\infty_{n=1}Z_n\delta_{(X_n,Y_n)}
\overset{d}{=}\sum^\infty_{n=1}Z'_n\delta_{(X'_n,Y'_n)}=:\zeta',
\end{align}
where $\sum^\infty_{n=1}\delta_{Z'_n}$ has a Poisson-Dirichlet distribution
with parameter $\alpha$ and $((X'_n,Y'_n))_{n\ge 1}$ is an independent
i.i.d.\ sequence with marginal distribution $\nu'\otimes\lambda$. 

Assume now that $\BX$ is a Borel space. Then, outside a measurable
$\BP$-null set, $\sum^\infty_{n=1}\delta_{Z_n}$ is a measurable function
of $\zeta$, while $\sum^\infty_{n=1}\delta_{Z'_n}$ is the result of the same 
measurable function applied to $\zeta'$; see e.g.\ \cite[Exercise 13.10]{LastPenrose17}.
(Here we use that $\BP(Y_m\ne Y_n)=1$ for $m\ne n$.)
Therefore \eqref{e5.22} shows that $\sum^\infty_{n=1}\delta_{Z_n}$ has 
the Poisson-Dirichlet distribution.\qed

\bigskip

If $\BX$ is a singleton, then Theorem \ref{t5.2} (and its converse)
reduces to a classical result from \cite{McCloskey65}, mentioned in the introduction.

\end{document}